\documentclass[a4paper,oneside]{amsart}

\usepackage{amsthm,amsfonts,amsmath,url}

\usepackage{fouriernc}

\newcommand{\dd}{\mathrm{d}}
\newcommand{\RR}{\mathbb{R}}
\newcommand{\ZZ}{\mathbb{Z}}
\newcommand{\TT}{\mathbb{T}}
\newcommand{\mx}{\mathrm{max}}

\theoremstyle{definition}

\begin{document}

\title[]{An inequality for the maximum curvature through a geometric flow}

\author{Konstantin Pankrashkin}

\address{Laboratoire de math\'ematiques (UMR 8628 du CNRS), Universit\'e Paris-Sud, B\^atiment 425, 91405 Orsay Cedex, France}

\email{konstantin.pankrashkin@math.u-psud.fr}


\begin{abstract} 
We provide a new proof of the following inequality: the maximum curvature $k_\mathrm{max}$ and the enclosed area $A$
of a smooth Jordan curve satisfy $k_\mathrm{max}\ge \sqrt{\pi/A}$. The feature of our proof
is the use of the curve shortening flow.

\bigskip

\noindent
 \tiny
This is a preliminary version. The final version will appear in Archiv der Mathematik published by Birkh\"auser (\url{http://www.springer.com/birkhauser/mathematics/journal/13}).
\end{abstract}

\subjclass{53A04, 53C44, 35B50}

\keywords{Maximum curvature, curve shortening flow, maximum principle, isoperimetric inequality}

\maketitle

The aim of the present note is to give a new proof
for the following inequality: if $\gamma\subset\RR^2$ is a smooth Jordan curve, then
\begin{equation}
    \label{eq-kmax}
k_\mx\ge \sqrt{\pi/A},
\end{equation}
where $k_\mx$ is the maximum curvature and $A$ is the enclosed area,
and the equality holds iff $\gamma$ is a circle. We remark that here and later on we work with the signed curvature,
in particular, the curvature is non-negative iff the curve is convex. The inequality \eqref{eq-kmax} follows from a result by Pestov-Ionin on
inscribed disks~\cite{pi}:
If the curvature of a smooth Jordan curve does not exceed some positive $\kappa>0$, then the interior 
of the curve contains a disk of radius $1/\kappa$. Hence, the comparison of the areas gives 
$\pi/\kappa^2\le A$, and \eqref{eq-kmax} is obtained for $\kappa=\kappa_\mx$.
The original work \cite{pi} is hardly available, and a complete proof can be found e.g. in~\cite[Proposition~2.1]{ba}.
We are going to show that the inequality \eqref{eq-kmax}
can be alternatively deduced from the properties of the curve shortening flow \cite{GH,gray}.
The use of geometric flows for isoperimetric inequalities is a well established machinery,
see e.g.~\cite{topp,meanbook}, but the link to the inequality \eqref{eq-kmax} seems to be new.
We also mention that Eq.~\eqref{eq-kmax} plays a role for Faber-Krahn-type inequalites
for some eigenvalue problems~\cite{pp}.
Our proof naturally splits in several parts.

{\it A. Uniqueness.} We remark first that if the inequality \eqref{eq-kmax}
is proved, then one may show in a standard way that the equality holds only for the circles, see e.g. \cite[Proposition~7]{pp};
we include the argument for the sake of completeness. By contradiction, assume that one has the equality in \eqref{eq-kmax}
for some $\gamma$ different from a circle.
At some point of $\gamma$ the curvature is strictly smaller than $k_\mx$, and
by a small local deformation of $\gamma$ near such a point we may construct a new smooth Jordan curve $\gamma'$
having the same maximum curvature $k'_\mx=k_\mx$ but enclosing a strictly smaller area $A'$, which gives
$k'_\mx<\sqrt{\pi/A'}$ and contradicts \eqref{eq-kmax}.

{\it B. The inequality holds for the star-shaped curves.} It is elementary to show  \eqref{eq-kmax}
for star-shaped curves, cf. e.g.~\cite[Theorem 2]{pp}, and we include the proof for convenience.
Assume that $\gamma$ is star-shaped with respect to the origin and denote by $\ell$ its length.
Let $\Gamma:\RR/\ell\ZZ\to \gamma\subset \RR^2$ be a properly oriented arc-length parametrization,
then the Frenet formula $\Gamma''=-k n$, where $k$ is the curvature and $n$ is the outer unit normal,
and the integration by parts give
\[
\int_0^\ell k \Gamma\cdot n\, \dd s=
-\int_{0}^\ell \Gamma\cdot \Gamma'' \, \dd s
=\int_{0}^\ell | \Gamma'|^2  \dd s=\ell.
\]
As $\gamma$ is star-shaped, we have  $ \Gamma \cdot n\ge 0$ and
\[
\ell=\int_0^\ell k  \Gamma \cdot n\, \dd s
\le k_\mx\int_0^\ell \Gamma \cdot  n \, \dd s
=k_\mx\int_\gamma x_1\dd x_2-x_2\,\dd x_1=2k_\mx A,
\]
and \eqref{eq-kmax} follows from the classical isoperimetric inequality
$\ell^2\ge 4\pi A$, see e.g. \cite[\S 2.10]{BZ}.

{\it C. Some properties of the flow by curvature.} The study of general curves will be reduced to the star-shaped ones using the flow by curvature
(also called the curve shortening flow). Denote $\TT:=\RR/\ZZ$ and let
$C(\cdot,0):\TT\to \RR^2$ be a smooth embedded curve. By \cite[Main theorem and introduction]{gray}, there exist $T>0$ and $C:\TT\times [0,T)\to\RR^2$ such that, for any $t$,
$C(\cdot,t)$ is a smooth embedded curve and
\[
\dfrac{\partial C(x,t)}{\partial t}=-k(x,t) n(x,t),
\]
where $n(x,t)$ and $k(x,t)$ are respectively the outer unit normal 
and the curvature of the curve $C(\cdot,t)$ at the point $C(x,t)$,
and the limiting shape is a round point, with convergence in $C^\infty$ norm, and, in
particular, there exists $\tau\in [0,T)$ such that $C(\cdot,t)$ is convex for $t\in[\tau,T)$.
We remark that a compact proof  can be found  in \cite{crelle}.
The following properties will be used,
see Section~1 in\cite{gray}: the area $A(t)$ enclosed by the curve $C(\cdot,t)$ is
\begin{equation}
     \label{eq-aa}
A(t)=A(0)-2\pi t,
\end{equation}
hence, $T=A(0)/(2\pi)$, and the curvature satisfies
\begin{equation}
   \label{eq-kk}
\dfrac{\partial k}{\partial t}=\dfrac{\partial^2 k}{\partial s^2}+k^3,
\end{equation}
where $\partial/\partial s$ means the derivative with respect to the arc-length
on $C(\cdot,t)$. 

{\it D. Proof of~\eqref{eq-kmax} for general curves.}
By using a suitable scaling we may assume that $A=\pi$,
then the sought inequality becomes $k_\mx\ge 1$.
By contradiction, assume that for some  $\gamma$ we have
\begin{equation}
 \label{eq-contr}
k_\mx<1
\end{equation}
and construct a family $C(\cdot,t)$, $t\in[0,1/2)$,
of curves evolving by curvature as in the part C with $C(\cdot,0)=\gamma$.
By \eqref{eq-aa}, the curve $C(\cdot, t)$
encloses the area $\pi(1-2t)$, hence, the enlarged curves
\[
\Sigma(\cdot, t):= \dfrac{1}{\sqrt{1-2t}} \,C(\cdot,t)
\]
enclose the constant area $\pi$, and the curvature $K$ on $\Sigma (\cdot, t)$
is  $K=\sqrt{1-2t} \,k$. 
Using the equality \eqref{eq-kk} we arrive at
\begin{equation}
      \label{eq-m1}
\dfrac{\partial K}{\partial t}
= \sqrt{1-2t}\,\dfrac{\partial^2 k}{\partial s^2} -\dfrac{1}{1-2t}\, K (1 - K^2).
\end{equation}
By \eqref{eq-contr}, there is $M\in(0,1)$ with $K(x,0)<M$ for all $x\in\TT$.
Let us show that
\begin{equation}
     \label{eq-km}
  K(x,t)<M<1 \text{ for all } (x,t). 
\end{equation}
Assume by contradiction that the inequality \eqref{eq-km} is false, then there exists
a minimal value $t_*\in(0,1/2)$ for which one can find $x_*\in\TT$ with $K(x_*,t_*)=M$,
and then $x_*$ is a maximum of  $K(\cdot,t_*)$. As a consequence
it is also a maximum of $k(\cdot,t_*)$, in particular,
$\partial^2 k/\partial s^2(x_*,t_*)\le 0$,
and the equality~\eqref{eq-m1} gives
\[
\dfrac{\partial K}{\partial t}(x_*,t_*)\le -\dfrac{1}{1-2t_*}\, M (1 - M^2)<0.
\]
It follows that for small positive $\varepsilon$ one has
$K(x_*,t_*-\varepsilon)>M$, which contradicts the above choice of $t_*$.
Hence, the claim~\eqref{eq-km} holds.
On the other hand, by part C, for some $\tau>0$ the curve
$\Sigma(\cdot,\tau)$ is convex, and, by part B, for some $x\in\TT$ we have $K(x,\tau)\ge 1$,
which contradicts the inequality~\eqref{eq-km}.
Therefore, the condition~\eqref{eq-contr} cannot be satisfied.

\section*{Acknowledgments}
The research was partially supported by ANR NOSEVOL (2011 BS01019 01) and GDR CNRS 2279 DynQua.
The author thanks Yurii Nikonorov for bringing reference~\cite{pi} to the attention.

\end{document}